\documentclass[12pt]{amsart}

\usepackage{amsmath,xspace,amssymb,mathrsfs}
\usepackage{color}

\input xy
\xyoption{all}
\xyoption{2cell}
\UseAllTwocells
\CompileMatrices

\newcommand{\Spec}{\operatorname{Spec}}
\renewcommand{\phi}{\varphi}

\newcommand{\Ker}{\operatorname{Ker}}

\newcommand{\Min}{\operatorname{Min}}

\newcommand{\Max}{\operatorname{Max}}

\newcommand{\Ann}{\operatorname{Ann}}

\newcommand{\Con}{\operatorname{C}}
\newcommand{\Su}{\operatorname{S}}

\newtheorem{proposition}{Proposition}[section]
\newtheorem{lemma}[proposition]{Lemma} 
\newtheorem{corollary}[proposition]{Corollary}
\newtheorem{theorem}[proposition]{Theorem}

\newtheorem{prop-def}[proposition]{Proposition and definition}

\theoremstyle{definition}

\newtheorem{remark}[proposition]{Remark}

\begin{document}

\title[p.p. rings and their generalizations]{Structure theory of p.p. rings and their generalizations}

\author[A. Tarizadeh]{Abolfazl Tarizadeh}
\address{ Department of Mathematics, Faculty of Basic Sciences,
University of Maragheh \\
P. O. Box 55136-553, Maragheh, Iran.
 }
\email{ebulfez1978@gmail.com}


\begin{abstract} In this paper, new and significant advances on the understanding the structure of p.p. rings and their generalizations have been made. Especially among them, it is proved that a commutative ring $R$ is a generalized p.p. ring if and only if $R$ is a generalized p.f. ring and its minimal spectrum is Zariski compact, or equivalently, $R/\mathfrak{N}$ is a p.p. ring and $R_{\mathfrak{m}}$ is a primary ring for all $\mathfrak{m}\in\Max(R)$. Some of the major results of the literature either are improved or are proven by new methods. In particular, we give a new and quite elementary proof to the fact that a commutative ring $R$ is a p.p. ring if and only if $R[x]$ is a p.p. ring.
\end{abstract}

\maketitle

\section{Introduction and Preliminaries}

In this paper all rings are commutative. Recall that a ring $R$ is said to be a p.p. ring if every principal ideal of  $R$ is a projective $R$-module, (for the free module case see Proposition \ref{Proposition II}). It is well known that the annihilator of every finitely generated projective module is generated by an idempotent element, for its proof see e.g. \cite[Corollary 3.2]{A. Tarizadeh 5}. Using this, then a ring $R$ is a p.p. ring if and only if for each $f\in R$, then $\Ann(f)$ is generated by an idempotent element of $R$. \\

The p.f. ring notion is the natural generalization of p.p. ring notion. In fact, a ring $R$ is said to be a p.f. ring if every principal ideal of $R$ is a flat $R$-module. Every p.p. ring is a p.f. ring, because every projective module is flat. But there are p.f. rings which are not p.p. rings. For example $\Con(X)$, the ring of real-valued continuous functions on $X:=\beta\mathbb{R}^{+}\setminus\mathbb{R}^{+}$, is a p.f. ring which is not a p.p. ring, for the details see \cite{Al Ezeh 3}. \\

It can be shown that every p.f. ring is a reduced ring. In fact, p.f. rings and reduced mp-rings are the same, see e.g. \cite[Theorem 6.4]{Tarizadeh},
(remember that a ring $R$ is said to be an mp-ring if each prime ideal of $R$ contains a unique minimal prime ideal of $R$. This is equivalent to the statement that $\mathfrak{p}+\mathfrak{q}=R$ for every distinct minimal prime ideals $\mathfrak{p}$ and $\mathfrak{q}$ of $R$). So it is natural to investigate similar notions for general rings (not necessarily reduced). Hirano \cite{Hirano} generalized p.p. ring notion by defining that a ring $R$ is said to be a generalized p.p. ring (or, GPP-ring) if for each $f\in R$ there exists a natural number $n\geqslant1$ such that $Rf^{n}$ is $R$-projective. Motivated the work of  Hirano, in \cite{Al Ezeh 4}, the author defined a ring to be a generalized p.f. ring (or, GPF-ring) if for each $f\in R$ there exists some $n\geqslant1$ such that $Rf^{n}$ is a flat $R$-module. Every GPP-ring is a GPF-ring. The converse is not true. For example, let $R$ be a p.f. ring which is not a p.p. ring, then $R\times\mathbb{Z}_{4}$ is a GPF-ring which is not a GPP-ring. The notion of quasi p.f. ring is a natural generalization of the GPF-ring notion, we study these rings in Sections 4 and 5. \\

In this paper, we continue the studies of
\cite{Tarizadeh}-\cite{Vasconcelos}. As an outcome, substantial progresses in the understanding the structure of p.p. rings and their generalizations have been made. Theorems \ref{Theorem II}, \ref{Theorem VI} and \ref{Theorem XI} provide interesting characterizations for p.p. rings. The joint collaboration and interplay  between concepts from commutative algebra and topology, is what makes the results and their proofs of this paper particularly interesting. Among many major results, Theorems \ref{Theorem XI} and \ref{Theorem IX} are the culmination of methods and results of this paper. The main result of \cite[Theorem 1.2]{Jondrup} was already proved by using the machinery of commutative algebra. In what follows, we give a simple and elementary proof for this major result, see Theorem \ref{Theorem VI}. In the same vein, we also prove that a ring $R$ is a p.f. ring if and only if the polynomial ring $R[x]$ is a p.f. ring, see Theorem \ref{Theorem V}. We also introduce the new notion of strongly purified ring and then Theorem \ref{Proposition IV} is obtained. Theorems \ref{Theorem Asena balam}, \ref{theorem t-a quasi p.f.}, \ref{Theorem 92629456}, \ref{Theorem VIII} and Corollary \ref{Coro v nice result} are further main results of this paper. \\

Let $I$ be an ideal of a ring $R$. Then $I$ is called a pure ideal if the canonical ring map $R\rightarrow R/I$ is a flat ring map. It can be shown that $I$ is a pure ideal if and only if for each $f\in I$ there exists some $g\in I$ such that $f(1-g)=0$. For its proof see e.g. \cite[Corollary 2.4]{A. Tarizadeh 5}. \\

The nil-radical of a ring $R$ is denoted by $\mathfrak{N}(R)$ or simply by $\mathfrak{N}$ if there is no confusion. Sometimes it is also denoted by $\sqrt{0}$. If $\mathfrak{p}$ is a minimal prime ideal of a ring $R$ and $f\in\mathfrak{p}$, then by passing to the localization $R_{\mathfrak{p}}$, we may find some $g\in R\setminus\mathfrak{p}$ such that $fg$ is nilpotent. \\

It can be shown that a ring $R$ is zero dimensional if and only if for each $f\in R$ there exists a natural number $n\geqslant1$ such that $f^{n}(1-fg)=0$ for some $g\in R$, for its proof see e.g. \cite[Theorem 2.2(v)]{A. Tarizadeh 4}. \\

Remember that a ring $R$ is called an absolutely flat ring if
each $R$-module is $R$-flat. It is well known that a ring $R$ is an absolutely flat ring if and only if it is von-Neumann regular ring (i.e., each $f\in R$ can be written as $f=f^{2}g$ for some $g\in R$). For its proof see e.g. \cite[Theorem 2.1]{A. Tarizadeh Mathematica}. \\

Let $R$ be a ring. Then $Z(R)=\{f\in R: \Ann(f)\neq0\}$ is called the set of zero-divisors of $R$. The localization $T(R):=S^{-1}R$ with $S=R\setminus Z(R)$ is called the total ring of fractions of $R$. \\

If $\mathfrak{p}$ is a minimal prime ideal of a ring $R$, then $0\notin S=(R\setminus\mathfrak{p})\big(R\setminus Z(R)\big)$. Thus there exists a prime ideal $\mathfrak{q}$ of $R$ which does not meet $S$. It follows that $\mathfrak{p}=\mathfrak{q}\subseteq Z(R)$. Therefore
$\bigcup\limits_{\mathfrak{p}\in\Min(R)}\mathfrak{p}\subseteq Z(R)$. If moreover $R$ is a reduced ring (or even, lessened ring in the sense of \cite[Definition 3.9]{A. Tarizadeh 4}) then the equality holds. \\

Recall that a ring is said to be a primary ring if its zero ideal is a primary ideal. For any ring $R$, we have $\mathfrak{N}(R)\subseteq Z(R)$. It is easy to see that a ring $R$ is a primary ring if and only if every zero-divisor element of $R$ is nilpotent, i.e., $\mathfrak{N}(R)=Z(R)$. \\

This paragraph is a key point in understanding some proofs of this paper. First note that the term ``regular ideal'' may have different meanings in the literature. But in this paper, by a regular ideal we mean an ideal which is generated by a set of idempotents. Let $I$ be a regular ideal of a ring $R$. If $f\in I$ then there exists an idempotent $e\in I$ such that $f=fe$. In particular, let $e\in R$ be an idempotent, if $f\in Re$ then $f=fe$. Moreover, if $e$ and $e'$ are idempotents of a ring $R$ such that $Re$ and $Re'$ are isomorphic as $R$-modules, then their annihilators are the same and so $e=e'$. Note that this does not hold in general. For example, $\mathbb{Z}\simeq2\mathbb{Z}$ but $1\neq2$. \\

If $\mathfrak{p}$ is a prime ideal of a ring $R$, then the canonical ring map $R\rightarrow R_{\mathfrak{p}}$ is denoted by $\pi_{\mathfrak{p}}$. \\

We say that a morphism of rings $\phi:R\rightarrow R'$ lifts idempotents if whenever $e'\in R'$ is an idempotent, then there exists an idempotent $e\in R$ such that $\phi(e)=e'$. In this case, we also say that the idempotents of $R$ can be lifted along $R'$ (via $\phi$). \\

\section{First properties of p.p. rings and their generalizations}

We have improved the following result by adding (iii) and (iv) as new equivalents. The equivalency of the classical criterion (ii) is also proved by a new method.

\begin{theorem}\label{Theorem II} For a ring $R$ the following statements are equivalent. \\
$\mathbf{(i)}$ $R$ is a p.p. ring. \\
$\mathbf{(ii)}$ $R$ is a p.f. ring and $T(R)$ is an absolutely flat ring. \\
$\mathbf{(iii)}$ $T(R)$ is an absolutely flat ring and the idempotents of $R$ can be lifted along each localization of $R$. \\
$\mathbf{(iv)}$ $T(R)$ is an absolutely flat ring and the idempotents of $R$ can be lifted along $T(R)$.
\end{theorem}

{\bf Proof.} $\mathbf{(i)}\Rightarrow\mathbf{(ii)},\mathbf{(iii)}:$ If $f\in R$ then there is an idempotent $e\in R$ such that $\Ann(f)=Re$. But $f-e$ is a non zero-divisor of $R$, because if $(f-e)g=0$ for some $g\in R$, then $fg=ge=ge^{2}=fge=0$ thus $g\in\Ann(f)$ and so $g=ge=0$. Therefore we may write $f/1=(f/1)^{2}z$ with $z:=1/(f-e)$. Hence, $T(R)$ is absolutely flat. Let $x=f/s\in S^{-1}R$ be an idempotent where $S$ is a multiplicative subset of $R$. The $R$-module $Rf$ is canonically isomorphic to $Re'$ where $e':=1-e$. Thus the ideal $(x)=(f/1)$ as $S^{-1}R$-module is canonically isomorphic to the ideal $(e'/1)$. Hence, $x=e'/1$. \\
$\mathbf{(ii)}\Rightarrow\mathbf{(i)}:$ If $f\in R$ then by hypothesis, $\Ann(f)$ is a pure ideal and there exists a non zero-divisor $s\in R$ such that $fs=f^{2}g$ for some $g\in R$. So there exists some $h\in\Ann(f)$ such that $s-fg=(s-fg)h$. Thus $(1-h)s=fg$. It follows that $h(1-h)s=0$ and so $h$ is an idempotent. If $h'\in\Ann(f)$ then $h'(1-h)s=0$ and so $h'=hh'$. Therefore $\Ann(f)=Rh$. \\ $\mathbf{(iii)}\Rightarrow\mathbf{(iv)}:$ There is nothing to prove. \\
$\mathbf{(iv)}\Rightarrow\mathbf{(i)}:$ If $f\in R$ then there exists some $x\in T(R)$ such that $f/1=(f^{2}/1)x$. Clearly $e':=(f/1)x$ is an idempotent and $\Ann_{T(R)}(f/1)=(1-e')$. By hypothesis, there is an idempotent $e\in R$ such that $e/1=1-e'$. It follows that $\Ann(f)=Re$. $\Box$ \\

\begin{proposition}\label{Proposition II} For a nonzero ring $R$ the following statements are equivalent. \\
$\mathbf{(i)}$ $R$ is an integral domain. \\
$\mathbf{(ii)}$ Every principal ideal of $R$ is a free $R$-module. \\
$\mathbf{(iii)}$ $R$ is a p.p. ring with trivial idempotents.
\end{proposition}

{\bf Proof.} $\mathbf{(i)}\Rightarrow\mathbf{(ii)}:$ Clearly $\Ann(f)=0$ for all nonzero $f\in R$, so $Rf\simeq R$. The zero ideal is also a free module with the empty basis. \\
$\mathbf{(ii)}\Rightarrow\mathbf{(i)}:$ If $f\in R$ then $Rf\simeq R/\Ann(f)$ is a free $R$-module. The annihilator of a free module is either the zero ideal or the whole ring. Thus $\Ann(f)=0$ or $\Ann(f)=R$. Hence, $R$ is an integral domain. \\
$\mathbf{(i)}\Rightarrow\mathbf{(iii)}:$ It is clear. \\
$\mathbf{(iii)}\Rightarrow\mathbf{(i)}:$ Suppose $fg=0$. If $f\neq0$ then $\Ann(g)=R$. So $g=0$. $\Box$ \\

\begin{lemma}\label{Lemma III} Let $R$ be a ring. If $f=\sum\limits_{i\geqslant0}r_{i}x^{i}\in R[[x]]$ is an idempotent, then $r_{0}$ is an idempotent and $r_{i}=0$ for all $i\geqslant1$.
\end{lemma}

{\bf Proof.} From $f=f^{2}$ we easily get that $r_{0}=r^{2}_{0}$. Let $k\geqslant1$ be the least natural number such that $r_{k}\neq0$. Then we will have $f=r_{0}+r_{k}x^{k}+\cdots$ and so $2r_{0}r_{k}=r_{k}$. This yields that $2r_{0}r_{k}=r_{0}r_{k}$ and so $r_{0}r_{k}=0$. We also have $r_{k}(1-r_{0})=0$. Thus $r_{k}=0$. $\Box$ \\

In the following result, we give a new and quite elementary proof to the main result of \cite[Theorem 1.2]{Jondrup}.

\begin{theorem}\label{Theorem VI} Let $R$ be a ring. Then $R$ is a p.p. ring if and only if $R[x]$ is a p.p. ring.
\end{theorem}

{\bf Proof.} Let $S:=R[x]$ be a p.p. ring. If $f\in R$ then by Lemma \ref{Lemma III}, there exists an idempotent $e\in R$ such that $\Ann_{S}(f)=Se$ and so $\Ann_{R}(f)=Re$. Conversely, let $R$ be a p.p. ring. If $f=\sum\limits_{i=0}^{n}f_{i}x^{i}\in S$ with the $f_{i}\in R$, then for each $i$ there exists an idempotent $e_{i}\in R$ such that $\Ann_{R}(f_{i})=Re_{i}$. It suffices to show that $\Ann_{S}(f)=Se$ where $e=\prod\limits_{i=0}^{n}e_{i}$. Clearly $Se\subseteq\Ann_{S}(f)$. Conversely, take $g=\sum\limits_{k=0}^{m}g_{k}x^{k}\in\Ann_{S}(f)$ with the $g_{k}\in R$. Then clearly $g_{0}=g_{0}e_{0}$. From $f_{0}g_{1}+f_{1}g_{0}=0$ we get that $f_{1}g_{0}e_{0}=0$ and so $g_{0}=g_{0}e_{0}e_{1}$. Then from $f_{2}g_{0}+f_{1}g_{1}+f_{0}g_{2}=0$ we obtain that $g_{0}=g_{0}e_{0}e_{1}e_{2}$. Thus by induction we will have $g_{0}=g_{0}e$. Therefore $f_{0}g_{1}=0$ and so $g_{1}=g_{1}e_{0}$ and by induction we obtain that $g_{1}=g_{1}e$. Thus $f_{0}g_{2}=0$ and by continuing this process we finally get that $g_{k}=g_{k}e$ for all $k=0,\ldots,m$. Hence, $g=ge$. $\Box$ \\

The following result also provides a simple proof to \cite[Remark on page 432]{Jondrup}.

\begin{theorem}\label{Theorem V} Let $R$ be a ring. Then $R$ is a p.f. ring if and only if $R[x]$ is a p.f. ring.
\end{theorem}

{\bf Proof.} Let $R$ be a p.f. ring. It is well known that a ring $R$ is a p.f. ring if and only if $R$ is a reduced mp-ring, see e.g. \cite[Theorem 6.4]{Tarizadeh}. Hence, $R$ is a reduced mp-ring and so $S:=R[x]$ is a reduced ring. If $P$ and $P'$ are distinct minimal prime ideals of $S$ then $\mathfrak{p}:=P\cap R$ and $\mathfrak{p}':=P'\cap R$ are distinct minimal prime ideals of $R$, since $P=\mathfrak{p}[x]$ and $P'=\mathfrak{p}'[x]$. We have $\mathfrak{p}+\mathfrak{p}'=R$ and so $P+P'=S$. Hence, $S$ is a reduced mp-ring and so it is a p.f. ring. Conversely, if $f\in R$ then it is easy to see that $\Ann_{R}(f)$ is a pure ideal. $\Box$ \\

The following result shows that one can easily construct generalized p.p. rings from a given ring.

\begin{proposition}\label{Corollary II} Let $\mathfrak{m}_{1},\ldots,\mathfrak{m}_{d}$ be finitely many distinct maximal ideals of a ring $R$ and $I=\bigcap\limits_{i=1}^{d}\mathfrak{m}^{c_{i}}_{i}$ with $c_{i}\geqslant1$. Then $R/I$ is a GPP-ring.
\end{proposition}

{\bf Proof.} If $\mathfrak{m}$ is a maximal ideal of a ring $R$ then $R/\mathfrak{m}^{k}$ is a GPP-ring for all $k\geqslant1$, because every element of the local ring $R/\mathfrak{m}^{k}$ is either invertible or nilpotent. It is also easy to see that the product of a finite family of rings is a GPP-ring if and only if each factor is a GPP-ring. Thus by the Chinese Remainder Theorem, $R/I\simeq R/\mathfrak{m}^{c_{1}}_{1}\times\cdots\times R/\mathfrak{m}^{c_{d}}_{d}$ is a GPP-ring. $\Box$ \\

\begin{corollary} If $n\in\mathbb{Z}$ then $\mathbb{Z}/n\mathbb{Z}$ is a GPP-ring.
\end{corollary}

{\bf Proof.} It follows from Proposition \ref{Corollary II}. $\Box$ \\

\begin{proposition} Let $R$ be a ring and $f\in R$. If $Rf^{n}$ is $R$-projective for some $n\geqslant1$, then $Rf^{k}$ as $R$-module is canonically isomorphic to $Rf^{n}$ for all $k\geqslant n$.
\end{proposition}

{\bf Proof.} It suffices to show that $\Ann(f^{k})=\Ann(f^{n})$ for all $k\geqslant n$. To see the latter it will be enough to show that $\Ann(f^{n})=\Ann(f^{n+1})$. There exists an idempotent $e\in R$ such that $\Ann(f^{n})=Re$. If $g\in\Ann(f^{n+1})$ then $fg=fge$. Thus $f^{n}g=f^{n-1}fg=0$ and so $g\in\Ann(f^{n})$. $\Box$ \\

\begin{proposition}\label{Proposition III} If $I$ is a pure ideal of a GPP-ring $R$, then $R/I$ is a GPP-ring.
\end{proposition}

{\bf Proof.} If $f\in R$ then there exists an idempotent $e\in R$ such that $\Ann(f^{n})=Re$ for some $n\geqslant1$. We show that $\Ann_{R/I}(f^{n}+I)=(Re+I)/I$. If $g+I\in\Ann_{R/I}(f^{n}+I)$ then there exists some $h\in I$ such that $(1-h)(1-e)g=0$. Thus $g\in Re+I$. $\Box$ \\

We provide an alternative proof to the following result which is the analogue of Proposition \ref{Proposition III}.

\begin{proposition}\cite[Theorem 1.11]{Al Ezeh 4}\label{Proposition V} If $I$ is a pure ideal of a GPF-ring $R$, then $R/I$ is a GPF-ring.
\end{proposition}

{\bf Proof.} If $f\in R$ then $\Ann(f^{n})$ is a pure ideal for some $n\geqslant1$. Thus $\big(\Ann(f^{n})+I\big)/I$ is a pure ideal, since the extension of a pure ideal under any ring map is a pure ideal. It suffices to show that $\Ann_{R/I}(f^{n}+I)=\big(\Ann(f^{n})+I\big)/I$. If  $g+I\in\Ann_{R/I}(f^{n}+I)$ then there exists some $h\in I$ such that $f^{n}g(1-h)=0$. Thus $g=g(1-h)+gh\in\Ann(f^{n})+I$. $\Box$ \\

\begin{lemma}\label{Lemma 7856 pure} Let $R$ be a ring and $f,g\in R$. If $\Ann(f)$ and $\Ann(g)$ are pure ideals, then $\Ann(fg)$ is a pure ideal.
\end{lemma}

{\bf Proof.} If $h\in\Ann(fg)$ then there exists some $f'\in\Ann(f)$ such that $(1-f')gh=0$. Thus there exists some $g'\in\Ann(g)$ such that $(1-f')(1-g')h=0$. So $(1-h')h=0$ where $h':=f'+g'-f'g'\in\Ann(fg)$. $\Box$ \\

\begin{corollary} Let $R$ be a ring and $f\in R$. If $Rf$ is $R-$flat, then $Rf^{n}$ is $R-$flat for all $n\geqslant0$.
\end{corollary}

{\bf Proof.} It follows from Lemma \ref{Lemma 7856 pure}. $\Box$ \\

\begin{corollary} The product of a finite family of rings is a GPF-ring if and only if each factor is a GPF-ring.
\end{corollary}

{\bf Proof.} It is deduced from Lemma \ref{Lemma 7856 pure}. $\Box$ \\

Let $R$ be either a local ring or an integral domain. Then clearly the zero ideal and the whole ring are the only pure ideals of $R$.

\begin{lemma}\label{Proposition Al-Ezeh} Every local GPF-ring is a primary ring.
\end{lemma}

{\bf Proof.} Let $R$ be a local GPF-ring. Take $f\in Z(R)$. It suffices to show that $f$ is nilpotent. By hypothesis, $\Ann(f^{n})$ is a pure ideal for some $n\geqslant1$. Since $R$ is a local ring, we have either $\Ann(f^{n})=0$ or $\Ann(f^{n})=R$.
But $\Ann(f)\subseteq\Ann(f^{n})$. Hence, we must have  $\Ann(f^{n})=R$. So $f^{n}=0$. $\Box$ \\

The only pure ideals of a primary ring are the zero ideal and the whole ring. In particular, every primary ring is a GPF-ring. \\

\begin{corollary}\label{Theorem VII} If $\mathfrak{p}$ is a prime ideal of a GPF-ring $R$, then $R_{\mathfrak{p}}$ is a primary ring.
\end{corollary}

{\bf Proof.} Clearly $R_{\mathfrak{p}}$ is a GPF-ring, because GPF-rings are stable under taking localizations. Thus by Lemma \ref{Proposition Al-Ezeh}, $R_{\mathfrak{p}}$ is a primary ring. $\Box$

In the following result, an alternative proof is given for \cite[Theorem 1.8]{Al Ezeh 4}.

\begin{theorem}\label{Theorem IV} For a ring $R$ the following statements are equivalent. \\
$\mathbf{(i)}$ $R$ is a GPF-ring. \\
$\mathbf{(ii)}$ If $f\in R$ then there exists a natural number $n\geqslant1$ such that in each localization $R_{\mathfrak{m}}$, either $f/1$ is a non zero-divisor or $f^{n}/1=0$.
\end{theorem}

{\bf Proof.} $\mathbf{(i)}\Rightarrow\mathbf{(ii)}:$ By hypothesis, there exists some $n\geqslant1$ such that $\Ann(f^{n})$ is a pure ideal. Let $\mathfrak{m}$ be a maximal ideal of $R$. If $\Ann(f^{n})$ is not contained in $\mathfrak{m}$, then in the ring $R_{\mathfrak{m}}$, $f^{n}/1=0$. If $\Ann(f^{n})\subseteq\mathfrak{m}$ then we show that in $R_{\mathfrak{m}}$, $f/1$ is a non zero-divisor. If not, then there exists a nonzero $g/1\in R_{\mathfrak{m}}$ such that $fgs=0$ for some $s\in R\setminus\mathfrak{m}$. This yields that $gs(1-h)=0$ for some $h\in\Ann(f^{n})$. So $s(1-h)\in\Ann(g)\subseteq\mathfrak{m}$ which is a contradiction. \\
$\mathbf{(ii)}\Rightarrow\mathbf{(i)}:$ We prove that $I=Rf^{n}$ is $R-$flat. It suffices to show that $I_{\mathfrak{m}}=(f^{n}/1)$ is
$R_{\mathfrak{m}}-$flat, since flatness is a local property. To see the latter it will be enough to show that
$J:=\Ann_{R_{\mathfrak{m}}}(f^{n}/1)$ is a pure ideal. If $f^{n}/1=0$ then $J=R_{\mathfrak{m}}$. If $f/1$ is a non zero-divisor, then $f^{n}/1$ is also a non zero-divisor and so $J=0$. $\Box$ \\

\section{p.p. rings and zero dimensionality}

In this section, the connections of p.p. rings and their generalizations with zero dimensional rings are studied. \\

Note that if $R$ is a reduced ring then for each $f\in R$, $\Ann(f)=\Ann(f^{n})$ for all $n\geqslant1$, because if $g\in\Ann(f^{n})$ then $fg$ is nilpotent and so $g\in\Ann(f)$. Hence, p.f. rings and reduced GPF-rings are the same. This easy argument proves
\cite[Theorem 1.9]{Al Ezeh 4}. The ring $\mathbb{Z}/4\mathbb{Z}$ is a GPF-ring which is not a p.f. ring, since it is not reduced (as another reason, $\Ann(2)=\{0,2\}$ is not a pure ideal).  \\

\begin{proposition} Every zero dimensional ring is a GPP-ring. In particular, every absolutely flat ring is a p.p. ring.
\end{proposition}

{\bf Proof.} Let $R$ be a zero dimensional ring. If $f\in R$ then it is well known that there exists a natural number $n\geqslant1$ such that $f^{n}(1-fg)=0$ for some $g\in R$, see e.g. \cite[Theorem 2.2(v)]{A. Tarizadeh 4}. It follows that $f^{n}=f^{2n}h$ for some $h\in R$. Thus $f^{n}h$ is an idempotent and $\Ann_{R}(f^{n})=R(1-f^{n}h)$. Hence, $R$ is a GPP-ring. It is also well know that a ring is absolutely flat if and only if it is reduced and zero dimensional. Clearly p.p. rings and reduced GPP-rings are the same. This argument establishes the second part of the assertion. $\Box$ \\

The following result  improves \cite[Theorem 1]{Hirano}.

\begin{lemma}\label{Theorem I} For a ring $R$ the following statements are equivalent.  \\
$\mathbf{(i)}$ $T(R)$ is zero dimensional. \\
$\mathbf{(ii)}$ If $f\in R$ then $\Ann(f^{n})=\Ann(f^{n+1})$ for some $n\geqslant1$, and there exists some $h\in\Ann(f^{n})$ such that $\Ann(f^{n})\cap\Ann(h)=0$. \\
$\mathbf{(iii)}$ If $f\in R$ then there exist a natural number $n\geqslant1$ and a non zero-divisor $g\in R$ such that $f^{n}g=f^{2n}$.
\end{lemma}

{\bf Proof.} $\mathbf{(i)}\Rightarrow\mathbf{(ii)}:$ There exist a natural number $n\geqslant1$ and a non zero-divisor $s\in R$ such that $f^{n}s=f^{2n}g$ for some $g\in R$. If $g'\in\Ann(f^{n+1})$ then $f^{n}g's=(f^{n+1}g')(f^{n-1}g)=0$ and so $g'\in\Ann(f^{n})$. Setting $h:=s-f^{n}g$. Then clearly $h\in\Ann(f^{n})$, and if $h'\in\Ann(f^{n})\cap\Ann(h)$ then $h's=h'f^{n}g=0$ and so $h'=0$. \\
$\mathbf{(ii)}\Rightarrow\mathbf{(iii)}:$ Setting $g:=f^{n}-h$. If $gg'=0$ then $f^{2n}g'=f^{n}hg'=0$ and so $g'\in\Ann(f^{2n})=\Ann(f^{n})$. Thus $g'h=0$. Hence, $g'\in \Ann(f^{n})\cap\Ann(h)=0$. So $g$ is a non zero-divisor of $R$. \\
$\mathbf{(iii)}\Rightarrow\mathbf{(i)}:$ Straightforward. $\Box$ \\

Note that if in a ring $R$ we have $\Ann(f^{n})=\Ann(f^{n+1})$ for some $f\in R$, then $\Ann(f^{n})=\Ann(f^{k})$ for all $k\geqslant n$. \\

The following result is the corresponding of Theorem \ref{Theorem II}, and it improves \cite[Theorem 2]{Hirano} and \cite[Theorem 2.1]{Al Ezeh 4}.

\begin{theorem}\label{Theorem III} For a ring $R$ the following statements are equivalent.  \\
$\mathbf{(i)}$ $R$ is a GPP-ring. \\
$\mathbf{(ii)}$ $R$ is a GPF-ring and $T(R)$ is zero dimensional. \\
$\mathbf{(iii)}$ $T(R)$ is zero dimensional and the idempotents of $R$ can be lifted along each localization of $R$. \\
$\mathbf{(iv)}$ $T(R)$ is zero dimensional and the idempotents of $R$ can be lifted along $T(R)$.
\end{theorem}

{\bf Proof.} $\mathbf{(i)}\Rightarrow\mathbf{(ii)}:$ If $f\in R$ then by  hypothesis, there exist a natural number $n\geqslant1$ and an idempotent $e\in R$ such that $\Ann(f^{n})=Re$. So $f^{n}-e$ is a non zero-divisor of $R$. Thus by Lemma \ref{Theorem I}, $T(R)$ is zero dimensional. \\
$\mathbf{(ii)}\Rightarrow\mathbf{(i)}:$ If $f\in R$ then by Lemma \ref{Theorem I}, $\Ann(f^{n})=\Ann(f^{n+1})$ for some $n\geqslant1$, and there exists some $h\in\Ann(f^{n})$ such that $\Ann(f^{n})\cap\Ann(h)=0$. First we show that $\Ann(f^{n})$ is a pure ideal. By hypothesis, there exists a natural number $k\geqslant1$ such that $\Ann(f^{k})$ is a pure ideal. If $k\geqslant n$ then we are done, since $\Ann(f^{n})=\Ann(f^{k})$. If $k<n$ then $n\leqslant kd$ for some positive natural number $d$. So by Lemma \ref{Lemma 7856 pure}, $\Ann(f^{kd})=\Ann(f^{n})$ is a pure ideal. Thus there exists some $g\in\Ann(f^{n})$ such that $1-g\in\Ann(h)$. Then $g(1-g)\in\Ann(f^{n})\cap\Ann(h)=0$. So $g$ is an idempotent. If $g'\in\Ann(f^{n})$ then $g'(1-g)\in\Ann(f^{n})\cap\Ann(h)=0$. Therefore $\Ann(f^{n})=Rg$. \\
$\mathbf{(i)}\Rightarrow\mathbf{(iii)}:$ It is proven exactly like the implication (i)$\Rightarrow$(iii) of Theorem \ref{Theorem II}. \\ $\mathbf{(iii)}\Rightarrow\mathbf{(iv)}:$ There is nothing to prove. \\ $\mathbf{(iv)}\Rightarrow\mathbf{(i)}:$ It is proved exactly like the implication (iv)$\Rightarrow$(i) of Theorem \ref{Theorem II}. $\Box$ \\

Let $R$ be a ring. Then there exists a unique topology over $\Spec(R)$, known as the flat topology, such that the collection of $V(f)=\{\mathfrak{p}\in\Spec(R): f\in\mathfrak{p}\}$ with $f\in R$ forms a subbase for its opens. For more information we refer the interested reader to \cite{Tarizadeh 2} or \cite{A. Tarizadeh 6 new method}. \\

\begin{theorem}\label{Theorem Asena balam} Let $R$ be a ring. If $T(R)$ is zero dimensional, then $\Min(R)$ is Zariski compact.
\end{theorem}

{\bf Proof.} By \cite[Theorem 3.3]{Tarizadeh 3}, it will be enough to show that if $f\in R$ then $\Min(R)\cap D(f)$ is a flat open of $\Min(R)$. There exist a natural number $n\geqslant1$ and a non zero-divisor $g\in R$ such that $f^{n}(g-fh)=0$ for some $h\in R$. Then we show that $\Min(R)\cap D(f)=\Min(R)\cap V(g-fh)$. The inclusion $\Min(R)\cap D(f)\subseteq\Min(R)\cap V(g-fh)$ is clear from the above relation. To see the reverse inclusion, take $\mathfrak{p}\in\Min(R)\cap V(g-fh)$. If $f\in\mathfrak{p}$ then $g\in\mathfrak{p}$. But this is a contradiction, because every minimal prime is contained in the set of zero-divisors. $\Box$ \\

\begin{corollary}\label{Lemma II} If $R$ is a GPP-ring, then $\Min(R)$ is Zariski compact.
\end{corollary}

{\bf Proof.} By Theorem \ref{Theorem III}, the ring $T(R)$ is zero dimensional. Then apply Theorem \ref{Theorem Asena balam}. $\Box$ \\

\begin{corollary} Let $R$ be a ring. If $T(R[x])$ is zero dimensional, then $\Min(R)$ is Zariski compact.
\end{corollary}

{\bf Proof.} For any ring $R$, the minimal prime ideals of $R[x]$ are precisely of the form $\mathfrak{p}[x]$ where $\mathfrak{p}$ is a minimal prime ideal of $R$. Hence, the map $\phi:\Min(R)\rightarrow\Min(R[x])$ given by $\mathfrak{p}\rightsquigarrow\mathfrak{p}[x]$ is bijective.
This map is continuous, because if $f=\sum\limits_{i=0}^{n}f_{i}x^{i}\in R[x]$ with the $f_{i}\in R$, then
$\phi^{-1}(U)=\bigcup\limits_{i=0}^{n}U_{i}$ where $U=\Min(R[x])\cap D(f)$ and $U_{i}=\Min(R)\cap D(f_{i})$. The converse of $\phi$ is also continuous, because it is induced by the ring extension $R\subseteq R[x]$. Therefore $\phi$ is a homeomorphism.
By Theorem \ref{Theorem Asena balam}, $\Min(R[x])$ is Zariski compact. Hence, $\Min(R)$ is Zariski compact. $\Box$

\section{Quasi p.f. rings}

Recall from \cite[\S2]{A. Tarizadeh 5} that an ideal $I$ of a ring $R$ is called a quasi-pure ideal (or, N-pure ideal in the sense of \cite{Aghajani}) if for each $f\in I$ there exists some $g\in I$ such that $f(1-g)$ is nilpotent. By a quasi p.f. ring (or, mid ring in the sense of \cite{Aghajani}) we mean a ring $R$ such that $\Ann(f)$ is a quasi-pure ideal for all $f\in R$. In the following result, we give various characterizations for quasi p.f. rings, (this result also can be found in \cite{Aghajani} which has been obtained independently). In addition to this result, our recent work \cite[Theorem 3.16]{A. Tarizadeh 4} provides further characterizations for quasi p.f. rings.

\begin{theorem}\label{theorem t-a quasi p.f.} For a ring $R$ the following statements are equivalent. \\
$\mathbf{(i)}$ $R$ is a quasi p.f. ring. \\
$\mathbf{(ii)}$ $R_{\mathfrak{p}}$ is a primary ring for all $\mathfrak{p}\in\Spec(R)$. \\
$\mathbf{(iii)}$ $R_{\mathfrak{m}}$ is a primary ring for all $\mathfrak{m}\in\Max(R)$.\\
$\mathbf{(iv)}$ $\Ker\pi_{\mathfrak{p}}$ is a pure ideal for all $\mathfrak{p}\in\Min(R)$. \\
$\mathbf{(v)}$ $\Ker\pi_{\mathfrak{m}}$ is a primary ideal for all $\mathfrak{m}\in\Max(R)$.
\end{theorem}

{\bf Proof.} $\mathbf{(i)}\Rightarrow\mathbf{(ii)}:$ If $f/s$ is a zero-divisor element of $R_{\mathfrak{p}}$, then there is a nonzero element $g/1$ in $R_{\mathfrak{p}}$ such that $fgt=0$ for some $t\in R\setminus\mathfrak{p}$. By hypothesis, there is some $h\in\Ann(g)\subseteq\mathfrak{p}$ such that $ft(1-h)$ is nilpotent. Hence, its image $ft(1-h)/1$ under the canonical ring map $R\rightarrow R_{\mathfrak{p}}$ and so $f/s=\big(ft(1-h)/1\big)\big(1/st(1-h)\big)$ are nilpotent. \\
$\mathbf{(ii)}\Rightarrow\mathbf{(iii)}:$ There is nothing to prove. \\
$\mathbf{(iii)}\Rightarrow\mathbf{(i)}:$ Choose $f\in R$ and let $g\in\Ann(f)$. It suffices to show that $\Ann(f)+\Ann(g^{n})=R$ for some $n\geqslant1$, because in this case there will be some $h\in\Ann(f)$ such that $g^{n}(1-h)=0$, this shows that $g(1-h)$ is nilpotent. Suppose $\Ann(f)+\Ann(g^{n})\neq R$ for all $n\geqslant1$. Then $I=\Ann(f)+\sum\limits_{n\geqslant1}\Ann(g^{n})$ will be a proper ideal of $R$. So $I\subseteq\mathfrak{m}$ for some $\mathfrak{m}\in\Max(R)$. It follows that in $R_{\mathfrak{m}}$, $f/1\neq0$ and $g/1$ is not nilpotent. But we have $fg=0$ and so in $R_{\mathfrak{m}}$, $(f/1)(g/1)=0$. By hypothesis, $R_{\mathfrak{m}}$ is a primary ring. Hence, we reach a contradiction. \\
$\mathbf{(i)}\Rightarrow\mathbf{(iv)}:$ If $f\in\Ker\pi_{\mathfrak{p}}$ then $fg=0$ for some $g\in R\setminus\mathfrak{p}$. Then by hypothesis, $g(1-h')$ is nilpotent for some $h'\in\Ann(f)$. It follows that $g^{n}(1-h)=0$ for some $n\geqslant1$ and some $h\in\Ann(f)$. So $1-h\in\Ker\pi_{\mathfrak{p}}$ and $f\big(1-(1-h)\big)=fh=0$. \\
$\mathbf{(iv)}\Rightarrow\mathbf{(iii)}:$ It suffices to show that each zero-divisor $f/s\in R_{\mathfrak{m}}$ is nilpotent. If not, then there exists a minimal prime ideal $\mathfrak{p}$ of $R$ with $\mathfrak{p}\subseteq\mathfrak{m}$ such that $f/s\notin\mathfrak{p}R_{\mathfrak{m}}$. So $f\notin\mathfrak{p}$. Since $f/s$ is a zero-divisor, then there exists a nonzero $g/1\in R_{\mathfrak{m}}$ such that $fgt=0$ for some $t\in R\setminus\mathfrak{m}$. It follows that $g\in\Ker\pi_{\mathfrak{p}}$. Then by hypothesis, $g(1-h)=0$ for some $h\in\Ker\pi_{\mathfrak{p}}\subseteq\mathfrak{p}\subseteq\mathfrak{m}$. So $1-h\in\Ann(g)\subseteq\mathfrak{m}$. Hence, $1=(1-h)+h\in\mathfrak{m}$ which is a contradiction. \\
$\mathbf{(iv)}\Rightarrow\mathbf{(v)}:$ There exists a minimal prime $\mathfrak{p}$ of $R$ such that $\mathfrak{p}\subseteq\mathfrak{m}$. Thus by \cite[Theorem 3.16]{A. Tarizadeh 4}, $\Ker\pi_{\mathfrak{m}}=\Ker\pi_{\mathfrak{p}}$. But for any ring $R$, if $\mathfrak{p}$ is a minimal prime ideal of it, then by \cite[Lemma 3.1]{A. Tarizadeh 4}, $\Ker\pi_{\mathfrak{p}}$ is a primary ideal. \\
$\mathbf{(v)}\Rightarrow\mathbf{(iv)}:$ By \cite[Theorem 3.16]{A. Tarizadeh 4}, it suffices to show that if $\mathfrak{p}$ is a minimal prime ideal of $R$ and $\mathfrak{m}$ is a maximal ideal of $R$ with $\mathfrak{p}\subseteq\mathfrak{m}$, then $\Ker\pi_{\mathfrak{m}}=\Ker\pi_{\mathfrak{p}}$. The inclusion $\Ker\pi_{\mathfrak{m}}\subseteq\Ker\pi_{\mathfrak{p}}$ is clear. To see the reverse inclusion, take $f\in\Ker\pi_{\mathfrak{p}}$. Then $fg=0$ for some $g\in R\setminus\mathfrak{p}$. If $f\notin\Ker\pi_{\mathfrak{m}}$, then by hypothesis, $g\in\sqrt{\Ker\pi_{\mathfrak{m}}}\subseteq
\sqrt{\Ker\pi_{\mathfrak{p}}}\subseteq\mathfrak{p}$ which is a contradiction. $\Box$ \\

\begin{theorem}\label{Theorem 92629456} The following statements  hold. \\
$\mathbf{(i)}$ Every GPF-ring is a quasi p.f. ring. \\
$\mathbf{(ii)}$ Every quasi p.f. ring is a mp-ring.
\end{theorem}

{\bf Proof.} $\mathbf{(i)}:$ It follows from Theorem \ref{theorem t-a quasi p.f.} and Corollary \ref{Theorem VII}. \\
$\mathbf{(ii)}:$ Let $\mathfrak{p}$ and $\mathfrak{q}$ be distinct minimal prime ideals of a ring $R$. We have $0\in(R\setminus\mathfrak{p})(R\setminus\mathfrak{q})$. Thus there are $f\in R\setminus\mathfrak{p}$ and $g\in R\setminus\mathfrak{q}$ such that $fg=0$. If $R$ is a quasi p.f. ring, then there exists some $h\in\Ann(f)$ such that $g(1-h)$ is nilpotent. It follows that $h\in\mathfrak{p}$ and $1-h\in\mathfrak{q}$. Hence, $\mathfrak{p}+\mathfrak{q}=R$. $\Box$ \\

The above result provides examples of rings which are not GPF-rings. For instance, every local ring whose the number of minimal primes $\geqslant2$ is not a GPF-ring. \\

The following result generalizes \cite[Theorem 1.10]{Al Ezeh 4}.

\begin{corollary}\label{Corollary I} Let $R$ be a ring. Then $R$ is a mp-ring if and only if $R/\mathfrak{N}$ is a p.f. ring.
\end{corollary}

{\bf Proof.} By \cite[Theorem 6.4]{Tarizadeh}, p.f. rings and reduced mp-rings are the same. $\Box$ \\

We call a ring $R$ an \emph{admissible ring} if for each $f\in R$ then the set $X_{f}:=\{\mathfrak{m}\in\Max(R): f/1\in\mathfrak{N}(R_{\mathfrak{m}})\}$ is Zariski quasi-compact.

\begin{theorem}\label{Theorem VIII} If $R$ is an admissible and quasi p.f. ring, then it is a GPF-ring.
\end{theorem}

{\bf Proof.} It suffices to show that $R$ satisfies in the condition (ii)
of Theorem \ref{Theorem IV}. Fix $f\in R$. If $\mathfrak{m}\in X_{f}$ then there exist a natural number $d_{\mathfrak{m}}\geqslant1$ and some $s_{\mathfrak{m}}\in R\setminus\mathfrak{m}$ such that $s_{\mathfrak{m}}f^{d_{\mathfrak{m}}}=0$. Thus $X_{f}\subseteq\bigcup\limits_{\mathfrak{m}\in Y}D(s_{\mathfrak{m}})$. Since $R$ is an admissible ring, so there exist a natural number $N\geqslant1$ and finitely many elements $s_{1},...,s_{k}\in R$ such that $X_{f}\subseteq\bigcup\limits_{i=1}^{k}D(s_{i})$ and $s_{i}f^{N}=0$ for all $i$. By Theorem \ref{theorem t-a quasi p.f.}, $X_{f}=\{\mathfrak{m}\in\Max(R): f/1\in Z(R_{\mathfrak{m}})\}$. Therefore in each localization $R_{\mathfrak{m}}$, either $f/1$ is a non zero-divisor or $f^{N}/1=0$. $\Box$ \\

\section{Structural results on p.p. rings and their generalizations}

The following technical result provides another interesting characterization for p.p. rings. After proving this result we were informed that the literature, especially \cite[Theorem 4.2.10]{Glaz 3}, \cite[Theorem 2.11]{Glaz 2}, \cite[Proposition 2.6]{Matlis}, \cite[Proposition 3.4]{Vasconcelos}, contain proofs for this result (some of these proofs are incomplete, see the comments following \cite[Theorem 2.11]{Glaz 2}). Our proof uses only the standard techniques of commutative algebra and completely differs with the former proofs from the scratch.

\begin{theorem}\label{Theorem XI} Let $R$ be a ring. Then $R$ is a p.p. ring if and only if $R$ is a p.f. ring and $\Min(R)$ is Zariski compact.
\end{theorem}

{\bf Proof.} If $R$ is a p.p. ring then by Corollary \ref{Lemma II}, $\Min(R)$ is Zariski compact. Conversely, if $f\in R$ then $U:=\Min(R)\cap D(f)$ is Zariski clopen (both open and closed) subset of $\Min(R)$, because for any ring $R$ if $\mathfrak{p}\in E:=\Min(R)\cap V(f)$ then there exists some $h\in R\setminus\mathfrak{p}$ such that $fh$ is nilpotent and so $\mathfrak{p}\in\Min(R)\cap D(h)\subseteq E$. The ring $R$ is a (reduced) mp-ring and so we may consider the function $\gamma:\Spec(R)\rightarrow\Min(R)$ which sends each prime ideal $\mathfrak{p}$ of $R$ into $\gamma(\mathfrak{p})$, the unique minimal prime ideal of $R$ contained in $\mathfrak{p}$. By \cite[Theorem 6.3]{Tarizadeh} or by \cite[Theorem 3.26]{A. Tarizadeh 4}, the map $\gamma$ is Zariski continuous if and only if $\Min(R)$ is Zariski compact. Thus there exists an idempotent $e\in R$ such that $\gamma^{-1}(U)=V(e)$, because it is well known that for any ring $R$ then the map $r\rightsquigarrow V(r)=D(1-r)$ is a bijective function from the set of idempotents of $R$ onto the set of clopens of $\Spec(R)$. We show that $\Ann(f)=Re$. If $\mathfrak{p}$ is a minimal prime ideal of $R$ then we have either $f\in\mathfrak{p}$ or $e\in\mathfrak{p}$, because if $e\notin\mathfrak{p}$ then $\mathfrak{p}=\gamma(\mathfrak{p})\notin U$ and so $f\in\mathfrak{p}$. Thus $fe\in\mathfrak{p}$.
Hence,
$fe\in\bigcap\limits_{\mathfrak{p}\in\Min(R)}\mathfrak{p}=\sqrt{0}=0$. So $Re\subseteq\Ann(f)$. If $g\in\Ann(f)$ then $\gamma^{-1}(U')\subseteq V(1-e)$ where $U':=\Min(R)\cap D(g)$. Therefore $g(1-e)\in\bigcap\limits_{\mathfrak{p}\in\Min(R)}\mathfrak{p}=0$
and so $g=ge\in Re$.  $\Box$ \\

\begin{corollary}\label{Corollary III} Let $R$ be a ring. Then $R/\mathfrak{N}$ is a p.p. ring if and only if $R$ is a mp-ring and $\Min(R)$ is Zariski compact.
\end{corollary}

{\bf Proof.} If $R/\mathfrak{N}$ is a p.p. ring then it is reduced mp-ring and so $R$ is a mp-ring. Moreover by Corollary \ref{Lemma II}, $\Min(R)\simeq\Min(R/\mathfrak{N})$ is Zariski compact.
Conversely, by Corollary \ref{Corollary I}, $R/\mathfrak{N}$ is a p.f. ring. The space $\Min(R/\mathfrak{N})\simeq\Min(R)$ is Zariski compact. Thus by Theorem \ref{Theorem XI}, $R/\mathfrak{N}$ is a p.p. ring. $\Box$ \\

\begin{remark} Let $R$ be a GPP-ring. Then by Theorem \ref{Theorem 92629456} and Corollaries \ref{Lemma II} and \ref{Corollary III}, $R/\mathfrak{N}$ is a p.p. ring. Here we give a direct proof for it. If $f\in R$ then there exist a natural number $n\geqslant1$ and an idempotent $e\in R$ such that $\Ann(f^{n})=Re$. Obviously $J:=\Ann_{R/\mathfrak{N}}(f+\mathfrak{N})=
\Ann_{R/\mathfrak{N}}(f^{n}+\mathfrak{N})$ since $R/\mathfrak{N}$ is reduced. Thus it suffices to show that
$J=(Re+\mathfrak{N})/\mathfrak{N}$. Clearly
$(Re+\mathfrak{N})/\mathfrak{N}\subseteq J$. Conversely, take $g+\mathfrak{N}\in J$. We may write $g=g_{1}+g_{2}$ with $g_{1}\in Re$ and $g_{2}\in R(1-e)$. Thus $f^{n}g=f^{n}g_{2}\in\mathfrak{N}$ and so $g_{2}^{m}\in\Ann(f^{nm})=\Ann(f^{n})$ for some $m\geqslant1$. Therefore $g^{m}_{2}=g^{m}_{2}e=g^{m-1}_{2}g_{2}e=g^{m-1}_{2}g_{2}(1-e)e=0$. Hence, $g\in Re+\mathfrak{N}$. \\
\end{remark}

The following theorem is one of the main results of this paper. This result also considerably improves \cite[Proposition 4]{Hirano}.

\begin{theorem}\label{Theorem IX} For a ring $R$ the following statements are equivalent. \\
$\mathbf{(i)}$ $R$ is a GPP-ring. \\
$\mathbf{(ii)}$ $R$ is a GPF-ring and $\Min(R)$ is Zariski compact. \\
$\mathbf{(iii)}$ $R$ is a quasi p.f. ring and $\Min(R)$ is Zariski compact. \\
$\mathbf{(iv)}$ $R/\mathfrak{N}$ is a p.p. ring and $R_{\mathfrak{m}}$ is a primary ring for all $\mathfrak{m}\in\Max(R)$.
\end{theorem}

{\bf Proof.} $\mathbf{(i)}\Rightarrow\mathbf{(ii)}:$ By Corollary \ref{Lemma II}, $\Min(R)$ is Zariski compact. \\
$\mathbf{(ii)}\Rightarrow\mathbf{(i)}:$ If $f\in R$ then there exists a natural number $n\geqslant1$ such that $\Ann(f^{n})$ is a pure ideal. By Theorem \ref{Theorem 92629456}, $R$ is a mp-ring. Thus by Corollary \ref{Corollary III}, $R/\mathfrak{N}$ is a p.p. ring. So there exists an idempotent $e\in R$ such that $\Ann_{R/\mathfrak{N}}(f+\mathfrak{N})=(e+\mathfrak{N})$, since the idempotents of a ring can be lifted modulo its nil-radical. Hence, there exists a natural number $k\geqslant1$ such that $ef^{k}=0$. We may find a natural number $\ell\geqslant1$ such that $m:=n\ell\geqslant k$. Thus $Re\subseteq\Ann(f^{k})\subseteq\Ann(f^{m})$.
If $z\in\Ann(f^{m})$ then there exists some $g\in\Ann(f^{m})$ such that $z=zg$, because by Lemma \ref{Lemma 7856 pure}, $\Ann(f^{m})$ is a pure ideal. Clearly $fg$ and so
$g(1-e)$ are nilpotent. This yields that $g^{d}=g^{d}e$ for some $d\geqslant1$. Therefore $z=zg^{d}=zg^{d}e\in Re$. Hence, $\Ann(f^{m})=Re$. \\
$\mathbf{(ii)}\Rightarrow\mathbf{(iii)}:$ By Theorem \ref{Theorem 92629456}, $R$ is a quasi p.f. ring. \\
$\mathbf{(iii)}\Rightarrow\mathbf{(iv)}:$ By Theorem \ref{Theorem 92629456}, $R$ is an mp-ring. Thus by Corollary \ref{Corollary III}, $R/\mathfrak{N}$ is a p.p. ring. By Theorem \ref{theorem t-a quasi p.f.}, each $R_{\mathfrak{m}}$ is a primary ring. \\
$\mathbf{(iv)}\Rightarrow\mathbf{(ii)}:$ By Corollary \ref{Corollary III}, $\Min(R)$ is Zariski compact. To prove that $R$ is a GPF-ring, it suffices to show that $R$ satisfies in the condition (ii) of Theorem \ref{Theorem IV}. If $f\in R$ then there exists an idempotent $e\in R$ such that $\Ann_{R/\mathfrak{N}}(f+\mathfrak{N})=(e+\mathfrak{N})$, because the idempotents of any ring can be lifted modulo its nil-radical. Thus
 there exists a natural number $n\geqslant1$ such that $Re\subseteq\Ann(f^{k})\subseteq\sqrt{Re}$ for all $k\geqslant n$. Take a maximal ideal $\mathfrak{m}$ of $R$. If $e\notin\mathfrak{m}$ then clearly in $R_{\mathfrak{m}}$, $f^{n}/1=0$. Thus it remains to show that if $e\in\mathfrak{m}$ then $f/1$ is a non zero-divisor of $R_{\mathfrak{m}}$. If not, then $f/1$ will be nilpotent, because $R_{\mathfrak{m}}$ is a primary ring. So there exist some $s\in R\setminus\mathfrak{m}$ and a natural number $d\geqslant1$ such that $sf^{d}=0$. It follows that $s\in\sqrt{Re}\subseteq\mathfrak{m}$ which is a contradiction. Hence, $f/1$ is a non zero-divisor
 of $R_{\mathfrak{m}}$. $\Box$ \\


It is very important to notice that in Theorem \ref{Theorem IX}(iv) the ``primary ring'' assumption is crucial. In other words, if $R/\mathfrak{N}$ is a p.p. ring then $R$ is not necessarily a GPP-ring.
For example, take $R=\mathbb{Z}[x]/I$ where $I=(x^{2},2x)$. Clearly $\mathfrak{p}=(x)/I$ is the only minimal prime ideal of $R$. Hence, $R/\mathfrak{N}=R/\mathfrak{p}$ is a p.p. ring (in fact, $R/\mathfrak{N}\simeq\mathbb{Z}$). But $R$ is not a generalized p.p. ring, because take $f=2+I$ then $\Ann_{R}(f^{n})=\mathfrak{p}=R(x+I)$ for all $n\geqslant1$, it is easy to see that $R$ has no nontrivial idempotents, thus the proper and nonzero ideal $\Ann_{R}(f^{n})$ can not be generated by an idempotent element. This example also shows that the converse of Theorem \ref{Theorem 92629456}(ii) does not hold, because $R$ is clearly an mp-ring, but by Theorem \ref{Theorem IX}, it is not a quasi p.f. ring (and hence it is not a GPF-ring). We learned this nice example from Pierre Deligne. \\

\begin{corollary}\label{Coro v nice result} Let $R$ be a ring which has either finitely many maximal ideals or its minimal spectrum is Zariski compact. Then the following statements are equivalent. \\
$\mathbf{(i)}$ $R$ is a GPP-ring. \\
$\mathbf{(ii)}$ $R$ is a GPF-ring. \\
$\mathbf{(iii)}$ $R$ is a quasi p.f. ring.
\end{corollary}

{\bf Proof.} The implications $\mathbf{(i)}\Rightarrow\mathbf{(ii)}\Rightarrow\mathbf{(iii)}$ hold for any ring $R$. In fact, for the implication $\mathbf{(ii)}\Rightarrow\mathbf{(iii)}$ see Theorem \ref{Theorem 92629456}(i). \\
$\mathbf{(iii)}\Rightarrow\mathbf{(i)}:$ If $\Min(R)$ is Zariski compact then by Theorem \ref{Theorem IX}, $R$ is a GPP-ring. If $R$ has finitely many maximal ideals, then by Theorem \ref{Theorem VIII}, $R$ is a GPF-ring. So if $f\in R$ then $Rf^{n}$ is $R$-flat for some $n\geqslant1$. It is well known that if a ring $R$ has finitely many maximal ideals, then every finitely generated flat module over $R$ is projective, for its proof see e.g. \cite[Theorem 3.11]{A. Tarizadeh 8 f.g. proj}. Hence, $Rf^{n}$ is $R$-projective. $\Box$ \\

In \cite{A. Tarizadeh 7 minimal} the compactness of the minimal spectrum with respect to the Zariski topology is characterized. \\

\begin{corollary}\label{coro iv quasi p.f.} Every local quasi p.f. ring is a primary ring.
\end{corollary}

{\bf Proof.} Let $R$ be a local quasi p.f. ring. Then by Corollary \ref{Coro v nice result}, $R$ is a GPF-ring. So by Lemma \ref{Proposition Al-Ezeh}, $R$ is a primary ring. $\Box$ \\

\begin{corollary} Every p.f. ring with finitely many maximal ideals is a p.p. ring.
\end{corollary}

{\bf Proof.} It follows from Corollary \ref{Coro v nice result}. $\Box$ \\

The following result shows that the quasi p.f. ring notion is a local property.

\begin{corollary} For a ring $R$ the following statements are equivalent. \\
$\mathbf{(i)}$ $R$ is a quasi p.f. ring. \\
$\mathbf{(ii)}$ $R_{\mathfrak{p}}$ is a quasi p.f. ring for all $\mathfrak{p}\in\Spec(R)$. \\
$\mathbf{(iii)}$ $R_{\mathfrak{m}}$ is a quasi p.f. ring for all $\mathfrak{m}\in\Max(R)$.
\end{corollary}

{\bf Proof.} $\mathbf{(i)}\Rightarrow\mathbf{(ii)}:$ Straightforward. In fact, if $S$ is a multiplicative subset of a quasi p.f. ring $R$, then $S^{-1}R$ is a quasi p.f. ring. \\
$\mathbf{(ii)}\Rightarrow\mathbf{(iii)}:$ There is nothing to prove. \\
$\mathbf{(iii)}\Rightarrow\mathbf{(i)}:$ By Corollary \ref{coro iv quasi p.f.}, each $R_{\mathfrak{m}}$ is a primary ring. Thus by Theorem \ref{theorem t-a quasi p.f.}, $R$ is a quasi p.f. ring. $\Box$ \\

It is easy to see that p.p. rings, p.f. rings, GPP-rings and GPF-rings are stable under taking localizations. \\

\begin{remark} Note that Theorems \ref{Theorem Asena balam}, \ref{Theorem XI} and \ref{Theorem IX} provide alternative proofs for
the implication ``(ii)$\Rightarrow$(i)'' of Theorems \ref{Theorem II} and \ref{Theorem III}. \\
\end{remark}

There is another natural and interesting generalization of p.p. ring notion which is called almost p.p. ring. In fact, a ring $R$ is said to be an almost p.p. ring if for each $f\in R$, then $\Ann(f)$ is a regular ideal. Obviously every p.p. ring is an almost p.p. ring, but $\Con(\beta\mathbb{N}\setminus\mathbb{N})$ is an example of  almost p.p. ring which is not a p.p. ring, for the details see \cite{Al Ezeh 3}. \\

In \cite[\S8]{Tarizadeh} we introduce and study purified rings. In fact, a ring $R$ is called a purified ring if for every distinct minimal prime ideals $\mathfrak{p}$ and $\mathfrak{q}$ of $R$, then there exists an idempotent $e\in\mathfrak{p}$ such that $1-e\in\mathfrak{q}$. Especially in \cite[Theorem 8.5]{Tarizadeh}, various characterizations for reduced purified rings are given. In particular, ``reduced purified ring'' and ``almost p.p. ring'' are the same. Using this and \cite[Theorem 8.3]{Tarizadeh}, then a ring $R$ is a purified ring if and only if $R/\mathfrak{N}$ is an almost p.p. ring. \\

It is also well known that a ring $R$ is almost p.p. ring if and only if $R[x]$ is almost p.p. ring, for its proof see \cite[Corollary 4.7]{Kosan et al}. It also easily follows from \cite[Theorem 8.5]{Tarizadeh}.
Note that in the commutative case,  ``feebly Baer ring'' and ``almost p.p. ring'' are the same. \\

Now we call a ring $R$ \emph{strongly purified} if for each $f\in R$ there exists a natural number $n\geqslant1$ such that $\Ann(f^{n})$ is a regular ideal. This definition generalizes at once both ``GPP-ring'' and ``almost p.p. ring'' notions. Furthermore, every strongly purified ring is a GPF-ring. \\

\begin{theorem}\label{Proposition IV} If $R$ is a strongly purified ring, then $R$ is a purified ring and every pure ideal of $R$ is a regular ideal.
\end{theorem}

{\bf Proof.} If $\mathfrak{p}$ and $\mathfrak{q}$ are distinct minimal prime ideals of $R$ then there exist $f\in R\setminus\mathfrak{p}$ and $g\in R\setminus\mathfrak{q}$
such that $fg=0$. By hypothesis, there exists a natural number $n\geqslant1$ such that $\Ann(f^{n})$ is generated by a set of idempotents of $R$. Thus there is an idempotent $e\in\Ann(f^{n})$ such that $g=ge$.
It follows that $e\in\mathfrak{p}$ and $1-e\in\mathfrak{q}$. Let $I$ be a pure ideal of $R$. If $f\in I$ then there exists some $g\in I$ such that $f(1-g)=0$. By hypothesis, there exists an idempotent $e\in R$ such that $f=fe$ and $e(1-g)^{n}=0$ for some $n\geqslant1$. It follows that $e\in I$. $\Box$ \\

\section{Ultra-products}

Let $(R_{x})$ be a family of rings indexed by a set $X$. If $f=(f_{x})\in R=\prod\limits_{x\in X}R_{x}$ then its support is defined as $\Su(f)=\{x\in X: f_{x}\neq0\}$. Let $I$ be an ideal of the power set ring $\mathcal{P}(X)$. Then it can be easily seen that $I^{\ast}=\{f\in R: \Su(f)\in I\}$ is an ideal of $R$. We call the quotient ring $R/I^{\ast}$ the ultra-ring (or, ultra-product) of the family $(R_{x})$ with respect to $I$. Ultra-rings, unlike direct products, possess many properties of the factor rings $R_{x}$. They are quite interesting and have vast applications in diverse fields of mathematics, see e.g. \cite{Barthel et al.}, \cite{Becker et al.}, \cite{Contessa}, \cite{Denef-Lipshitz}, \cite{Erman et al.} and \cite{Schoutens}. Our approach simplifies this construction in the literature. \\

\begin{lemma}\label{Lemma IV} Under the above hypothesis and terminology, then $I^{\ast}$ is a pure ideal of $R$.
\end{lemma}

{\bf Proof.} If $f\in I^{\ast}$ then consider the sequence $g=(g_{x})\in R$ where each $g_{x}$ is either $1$ or $0$, according as $x\in\Su(f)$ or $x\notin\Su(f)$. Then clearly $f=fg$ and
$\Su(g)\subseteq\Su(f)\in I$. So $g\in I^{\ast}$. $\Box$ \\

\begin{lemma}\label{Lemma V} If each $R_{x}$ is a reduced ring, then $R/I^{\ast}$ is a reduced ring.
\end{lemma}

{\bf Proof.} It is easy to see that $\Su(f)=\Su(f^{n})$ for all $f\in R$ and $n\geqslant2$. $\Box$ \\

\begin{proposition}\label{Theorem ultra xi} Let $(R_{x})$ be a family of rings indexed by a set $X$, let $R=\prod\limits_{x\in X}R_{x}$ and let $I$ be an ideal of $\mathcal{P}(X)$. Then the following statements hold. \\
$\mathbf{(i)}$ If each $R_{x}$ is a p.p. ring, then $R/I^{\ast}$ is a p.p. ring. \\
$\mathbf{(ii)}$ If each $R_{x}$ is a p.f. ring, then $R/I^{\ast}$ is a p.f. ring.
\end{proposition}

{\bf Proof.} $\mathbf{(i)}:$ Clearly $R$ is a p.p. ring, since it is easy to see that the direct product of a family of p.p. rings is a p.p. ring. Thus by Proposition \ref{Proposition III} and Lemma \ref{Lemma IV}, $R/I^{\ast}$ is a GPP-ring. But p.p. rings and reduced GPP-rings are the same. Thus by Lemma \ref{Lemma V}, $R/I^{\ast}$ is a p.p. ring. \\
$\mathbf{(ii)}:$ It is proven exactly like (i), by applying Proposition \ref{Proposition V} instead of Proposition \ref{Proposition III}. $\Box$ \\

\textbf{Acknowledgements.} The author would like to give sincere thanks to the referee for very careful reading of the paper. \\

\end{document}